\documentclass[11pt,reqno]{amsart}

\usepackage{cite,a4wide,amsmath,amssymb,graphicx,enumitem,amsaddr,tikz,pgfplots,hhline,epstopdf,multirow,multicol,soul}
\usepackage[hang,flushmargin]{footmisc}

\usetikzlibrary{decorations.markings}

\setlist[enumerate]{leftmargin=8mm}
\setlist[itemize]{leftmargin=8mm}

\newcommand{\cR}{\mathcal{R}}

\newcommand{\ds}{\frac{\text{d}S(t)}{\text{d}t}}
\newcommand{\di}{\frac{\text{d}I(t)}{\text{d}t}}
\newcommand{\dr}{\frac{\text{d}R(t)}{\text{d}t}}

\definecolor{darkgreen}{RGB}{34,177,76}

\DeclareMathOperator{\eff}{eff}

\begin{document}
\title{A design of governmental policies for the eradication of COVID-19 in Jakarta using an SIR-type mathematical model}

\author{\small Benny Yong, Jonathan Hoseana, Livia Owen}
\address{\normalfont\small Department of Mathematics, Parahyangan Catholic University, Bandung 40141, Indonesia}
\email{benny\_y@unpar.ac.id\textnormal{, }j.hoseana@unpar.ac.id\textnormal{, }livia.owen@unpar.ac.id}
\date{}

\begin{abstract}
Using a discretised version of our recently-developed SIR-type mathematical model for the spread of COVID-19, we construct a design of governmental policies for the eradication of the disease in the province of DKI Jakarta, Indonesia, taking as a basis the actual data of mid-2021. The design takes the form of a precise, quantitative method to determine the appropriate level(s) of restrictions on community activities (PPKM) which should be enforced in the province on any given day, based on the current values of the disease's effective reproduction number and the hospitals' bed-occupancy rate.

\smallskip\noindent\textsc{Keywords.} COVID-19; governmental policy; effective reproduction number; bed-occupancy rate

\smallskip\noindent\textsc{2020 MSC subject classification.} 92C60; 92D30
\end{abstract}

\maketitle

\section{Introduction}\label{section:Introduction}

Following the entry of the coronavirus disease 2019 (COVID-19) to Indonesia on March 1st, 2020 (infecting a 34-year-old woman after being in contact with a Japanese visitor who was tested positive only after leaving the country \cite{Ratcliffe,SetiawatyEtAl}), the disease spread across all 34 provinces \cite{WHOI16Apr2020}, the country's number of daily new cases reaching its currently highest point of 56,757 on July 15th, 2021 \cite{WHOI25Aug2021}. On July 3rd, the government called for \textit{emergency restrictions on community activities} (emergency PPKM) \cite{WHOI07Jul2021}, seemingly effectively creating a declining trend past the aforementioned highest point \cite{WHOI25Aug2021}. The emergency PPKM was subsequently referred to as a \textit{level 4} PPKM \cite{WHOI28Jul2021}, before, in the Special Capital Region of Jakarta (DKI Jakarta) province, downgraded to a \textit{level 3} PPKM on August 24th \cite{Bhwana}. Each PPKM level (1 to 4) signify a different degree of restrictions, the practical details being described in some instructional documents issued by the Minister of Home Affairs \cite{DepartemenDalamNegeri1,DepartemenDalamNegeri2}; for an (unofficial) English translation, see \cite{DoubleM}. Non-medical shopping centres, for instance, are instructed to shut completely on a level 4 PPKM, operate only at 25\% capacity with a 17:00 curfew on a level 3 PPKM, operate only at 50\% capacity with a 20:00 curfew on a level 2 PPKM, and operate only at 75\% capacity with a 21:00 curfew on a level 1 PPKM.

The government of every province, being the policy maker, is responsible for a judicious determination of the appropriate level of PPKM to be enforced on the respective province, on any given period. A prolonged high-level PPKM, albeit perhaps effective from the viewpoint of the disease's eradication, could lead to a significant decline in the economic sector. On the other hand, determining whether a low-level PPKM (or a level 0 PPKM, that is, a complete absence of restrictions) can already be allowed is extremely difficult.

In this paper, we formulate precise, quantitative conditions warranting the enforcement of each level of PPKM in the province of DKI Jakarta. These conditions are expressed in terms of ranges of values of the disease's effective reproduction number (the expected number of susceptible individuals to whom a single infected individual transmits the disease \cite{ChowellEtAl,CintronAriasEtAl,KeelingRohani}) and the hospitals' bed-occupancy rate (the number of hospitalised COVID-19 cases per isolation bed \cite{WHOI22Jul2020}), using the following mathematical model which we recently developed \cite{YongOwenHoseana}:
$$\left\{\begin{array}{rcl}
\displaystyle\ds &=& \displaystyle\lambda - \mu S(t) - \frac{\beta S(t)I(t)}{1+\gamma S(t)},\\[0.4cm]
\displaystyle\di &=& \displaystyle-\mu I(t) - \mu' I(t) +\frac{\beta S(t)I(t)}{1+\gamma S(t)}- \frac{\alpha I(t)}{1+\rho I(t)},\\[0.4cm]
\displaystyle\dr &=& \displaystyle-\mu R(t) + \frac{\alpha I(t)}{1+\rho I(t)}.
\end{array}\right.
$$
This is an SIR-type model, originally designed to capture the role of the susceptible individuals' cautiousness rate $\gamma\in[0,1]$ and the hospitals' bed-occupancy rate $\rho\in[0,1]$ in the disease's spread. The symbols $S(t)$, $I(t)$, $R(t)$ denote the number of susceptible, infected, recovered individuals at time $t\in[0,\infty)$, respectively, while $\lambda$, $\mu$, $\mu'$, $\beta$, $\alpha$ are positive parameters interpretable as the entry rate, the death rate, the death rate increment due to the disease, the incidence rate, and the recovery rate, respectively.

\section{Method and results}\label{sec:method}

\begin{figure}
\begin{tikzpicture} \pgfkeys{/pgf/number format/.cd,fixed,precision=0}
\begin{axis}[
	xmin=0,
	xmax=30,
	ymin=0.35,
	ymax=0.75,
    xtick={0,18,30},
    ytick={0.40,0.72},
    yticklabels={0.40,0.72},
	xlabel={$t$},
	ylabel={$\cR^{\eff}_t$},
	ylabel near ticks,
	samples=100,
	width=6cm,
	height=6cm,
	scaled y ticks=false
]
\addplot[only marks,mark=x,red] coordinates {(0,0.5365898354) (1,0.4214569248) (2,0.4015698535) (3,0.4319557268) (4,0.526801227) (5,0.6012665162) (6,0.6297432195) (7,0.5251848539) (8,0.4354241347) (9,0.4154420226) (10,0.4347198015) (11,0.4977222083) (12,0.6044268704) (13,0.6982226662) (14,0.6813821205) (15,0.6156235045) (16,0.5838196678) (17,0.5889674371) (18,0.6327882301) (19,0.6658684746) (20,0.6683216161) (21,0.6243536937) (22,0.613793243) (23,0.6260182539) (24,0.6725907329) (25,0.7172689424) (26,0.7159562328) (27,0.668974215) (28,0.5981084439) (29,0.5519767605) (30,0.5824943297)};

\addplot[only marks,mark=o,blue] coordinates {(0,0.5365898354) (1,0.4215353017) (2,0.4016814662) (3,0.4321456362) (4,0.5271471337) (5,0.6017370374) (6,0.6303222837) (7,0.5257503983) (8,0.4359539846) (9,0.4159764262) (10,0.4353158437) (11,0.498466812) (12,0.6053866704) (13,0.6994145166) (14,0.6826208182) (15,0.61679795) (16,0.584973232) (17,0.5901719339) (18,0.6341463316) (19,0.6673501948) (20,0.6698657468) (21,0.6258385862) (22,0.6152953039) (23,0.6275847834) (24,0.6743133452) (25,0.7191699339) (26,0.7178997226) (27,0.6708269693) (28,0.599798985) (29,0.5535650316) (30,0.5841945823)};

\draw[dashed] (axis cs:18,0.35) -- (axis cs:18,0.75);
\draw[dashed] (axis cs:0,0.4) -- (axis cs:30,0.4);
\draw[dashed] (axis cs:0,0.72) -- (axis cs:30,0.72);
\end{axis}
\end{tikzpicture}\qquad\begin{tikzpicture} \pgfkeys{/pgf/number format/.cd,fixed,precision=0}
\begin{axis}[
	xmin=0,
	xmax=30,
	ymin=0.1,
	ymax=0.45,
    xtick={0,18,30},
    ytick={0.1,0.2,0.3,0.4},
    yticklabels={0.1,0.2,0.3,0.4},
	xlabel={$t$},
	ylabel={$\rho_t$},
	ylabel near ticks,
	samples=100,
	width=6cm,
	height=6cm
]
\addplot[only marks,mark=x,red] coordinates {(0,0.43) (1,0.4) (2,0.39) (3,0.35) (4,0.34) (5,0.33) (6,0.32) (7,0.3) (8,0.27) (9,0.27) (10,0.26) (11,0.26) (12,0.24) (13,0.24) (14,0.23) (15,0.23) (16,0.22) (17,0.21) (18,0.2) (19,0.2) (20,0.18) (21,0.17) (22,0.18) (23,0.17) (24,0.17) (25,0.16) (26,0.15) (27,0.15) (28,0.15) (29,0.15) (30,0.13)};

\draw[dashed] (axis cs:18,0.1) -- (axis cs:18,0.45);
\draw[dashed] (axis cs:0,0.2) -- (axis cs:30,0.2);
\end{axis}
\end{tikzpicture}
\caption{\label{fig:Reff} On the left panel, a plot of the points $\left\{\left(t,\cR^{\eff}_t\right)\right\}_{t=0}^{30}$, where $\cR^{\eff}_t$ is computed using the formula \eqref{eq:Reff} (red crosses), and where $\cR^{\eff}_t$ is computed using the formula \eqref{eq:Reff} with $\hat S_t$ and $\hat N_t=\hat S_t+\hat I_t+\hat R_t$ replacing $S_t$ and $N_t$ (blue noughts). Here, $\hat\alpha_t$ and $\hat\beta_t$ are estimated using the L-BFGS-B algorithm. On the right panel, a plot of points $\left\{\left(t,\rho_t\right)\right\}_{t=0}^{30}$.}
\end{figure}
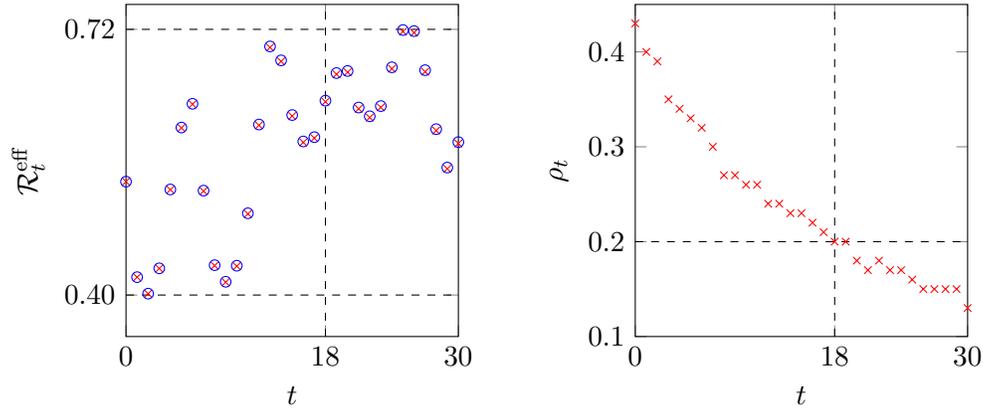

Let us now describe our method and results. Once again, our goal is to use the above model to formulate quantitative conditions under which each level of PPKM should be enforced in DKI Jakarta. To this end, we first discretise the above model as follows:
\begin{equation}\label{eq:discretemodel}
\left\{\begin{array}{rcl}
\hat S_{t+1} &=& \displaystyle \hat S_t+\lambda - \mu \hat S_t - \frac{\hat\beta_t \hat S_t \hat I_t}{1+\gamma \hat S_t},\\[0.4cm]
\hat I_{t+1} &=& \displaystyle \hat I_t-\mu \hat I_t - \mu' \hat I_t +\frac{\hat\beta_t \hat S_t \hat I_t}{1+\gamma \hat S_t}- \frac{\hat\alpha_t \hat I_t}{1+\rho_t \hat I_t},\\[0.4cm]
\hat R_{t+1} &=& \displaystyle \hat R_t-\mu \hat R_t + \frac{\hat\alpha_t \hat I_t}{1+\rho_t \hat I_t}.
\end{array}\right.
\end{equation}
We describe the quantities involved in or related to \eqref{eq:discretemodel}, and their assumed values.
\begin{enumerate}[topsep=5pt,itemsep=5pt]
\item[(i)] The following are parameters, whose values are therefore kept fixed:
\begin{equation}\label{eq:fixedparameters}
\lambda=\frac{10,467,629}{65\times 265}\approx 607.6997968,\,\,\,\mu=0.0000421496,\,\,\,\mu'=0.06,\,\,\,\text{and}\,\,\,\gamma=0.35,
\end{equation}
the first three of which are taken from \cite[Table 1]{AldilaEtAl}, while the last one is the estimate we used in \cite[Table 4]{YongOwenHoseana}.
\item[(ii)] Considering the available data, we let the independent variable $t\in\mathbb{Z}$ be measured in days, where $t=0$ represents August 7th, 2021. Thus, August 24th ---the day on which the PPKM was downgraded from level 4 to level 3--- is represented by $t=17$. Although we also consider negative values of $t$, we shall mainly work over the period $t\in\{0,\ldots,30\}$, i.e., from August 7th to September 6th.
\item[(iii)] From the data $\left\{\left(I_t,R_t\right)\right\}_{t=-90}^{30}$ of the daily numbers of infected and recovered individuals over the period from May 9th to September 6th, provided by the Johns Hopkins Coronavirus Resource Center \cite{JHU}, we obtain $\left\{\left(S_t,I_t,R_t\right)\right\}_{t=-90}^{30}$, where, for every $t\in\{-90,\ldots,30\}$, $S_t$ is the unique solution of the equation $S_t+I_t+R_t=N_t=10,467,629$, the right-hand side being the population size of DKI Jakarta, assumed to be constant and equal to the numerator of $\lambda$ in \eqref{eq:fixedparameters}. To form an initial condition $\left(\hat S_0,\hat I_0,\hat R_0\right)$ for \eqref{eq:discretemodel}, we let $\hat I_0=I_0=12,799$, $\hat R_0=R_0=802,334$, and $\hat S_0=S_0=N_0-I_0-R_0=9,652,496$.
\item[(iv)] The hospitals' bed-occupancy rate $\rho_t$ is assumed to be time-dependent; we use the data of daily occupancy percentages of hospital beds which are allocated for COVID-19 patients, provided by the Ministry of Health \cite{Kemkes}. See Figure \ref{fig:Reff} (right) for a plot of the points $\left\{\left(t,\rho_t\right)\right\}_{t=0}^{30}$.
\item[(v)] The data of the recovery and incidence rates, $\hat\alpha_t$ and $\hat\beta_t$, also assumed to be time-dependent, are unavailable. These quantities, therefore, must be estimated; see below. (This is the reason why we use the symbols $\hat\alpha_t$ and $\hat\beta_t$, rather than $\alpha_t$ and $\beta_t$.)
\end{enumerate}

As mentioned in the previous section, the conditions are to be formulated in terms of the model's effective reproduction number\footnote{We compute $\cR^{\eff}_t$ as $\cR_t\left(S_t/N_t\right)$ \cite{ChowellEtAl,CintronAriasEtAl,KeelingRohani}, where $\cR_t$ denotes the model's basic reproduction number \cite[equation (3)]{YongOwenHoseana}.} 
\begin{equation}\label{eq:Reff}
\cR^{\eff}_t=\frac{\hat\beta_t\lambda S_t}{(\mu+\gamma\lambda)\left(\mu+\mu'+\hat\alpha_t\right)N_t}
\end{equation}
and the hospitals' bed-occupancy rate $\rho_t$. Notice that $\cR^{\eff}_t$ and $\rho_t$ are independent, and thus complement each other. Indeed, we shall see later that, from the viewpoint of these two quantities, the downgrading of the PPKM level from 4 to 3 on August 24th \cite{Bhwana} is not at all due to a sustained decline of $\cR^{\eff}_t$, but to that of $\rho_t$.

We thus regard every point $\left(\cR^{\eff}_t,\rho_t\right)$ on the $\cR^{\eff}\rho$-plane as a quantitative description of the pandemic situation in DKI Jakarta on day $t$. Our aim is to assign to every such point the level(s) of PPKM (0, 1, 2, 3, or 4) which is appropriate to be enforced.

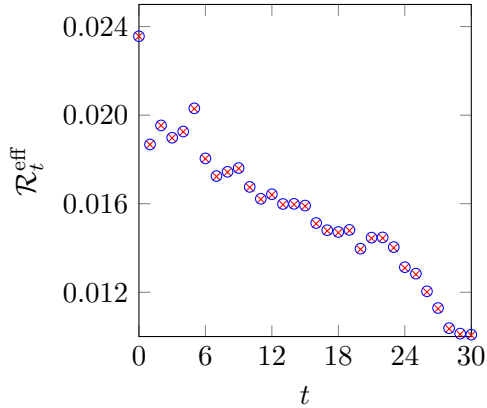
\begin{figure}
\begin{tikzpicture} \pgfkeys{/pgf/number format/.cd,fixed,precision=0}
\begin{axis}[
	xmin=0,
	xmax=30,
	ymin=0.01,
	ymax=0.025,
    xtick={0,6,12,18,24,30},
    ytick={0.012,0.016,0.020,0.024},
    yticklabels={0.012,0.016,0.020,0.024},
	xlabel={$t$},
	ylabel={$\cR^{\eff}_t$},
	ylabel near ticks,
	samples=100,
	width=6cm,
	height=6cm,
	scaled y ticks=false
]
\addplot[only marks,mark=x,red] coordinates {(0,0.023562588) (1,0.018672864) (2,0.019530949) (3,0.018970424) (4,0.019248796) (5,0.020290876) (6,0.018030524) (7,0.017228385) (8,0.017421162) (9,0.017587559) (10,0.01673497) (11,0.016197553) (12,0.016402392) (13,0.015960261) (14,0.015964694) (15,0.015886492) (16,0.015098998) (17,0.014773193) (18,0.014690805) (19,0.01478384) (20,0.013936279) (21,0.014434188) (22,0.014441651) (23,0.014004842) (24,0.013098541) (25,0.012806384) (26,0.012011179) (27,0.011257134) (28,0.010349349) (29,0.010102916) (30,0.010055584)};
\addplot[only marks,mark=o,blue] coordinates {(0,0.023562588) (1,0.018676336) (2,0.019536377) (3,0.018978765) (4,0.019261435) (5,0.020306755) (6,0.018047104) (7,0.017246937) (8,0.017442361) (9,0.017610183) (10,0.016757915) (11,0.016221785) (12,0.016428438) (13,0.015987505) (14,0.015993717) (15,0.0159168) (16,0.015128832) (17,0.014803406) (18,0.014722335) (19,0.014816737) (20,0.013968478) (21,0.014468516) (22,0.014476992) (23,0.014039887) (24,0.013132088) (25,0.012840325) (26,0.012043784) (27,0.011288311) (28,0.010378601) (29,0.010131986) (30,0.010084935)};
\end{axis}
\end{tikzpicture}
\caption{\label{fig:preliminaryReff} A plot of the points $\left\{\left(t,\cR^{\eff}_t\right)\right\}_{t=0}^{30}$, where $\cR^{\eff}_t$ is computed using the formula \eqref{eq:Reff} (red crosses), and where $\cR^{\eff}_t$ is computed using the formula \eqref{eq:Reff} with $\hat S_t$ and $\hat N_t=\hat S_t+\hat I_t+\hat R_t$ replacing $S_t$ and $N_t$ (blue noughts). Here, $\hat\alpha_t$ and $\hat\beta_t$ are estimated using \eqref{eq:estimates1}. The vertical distance between every two corresponding points is less than $10^{-4}$.}
\end{figure}

Let us now discuss the estimation of the recovery and incidence rates: $\hat\alpha_t$ and $\hat\beta_t$, over the period $t\in\{0,\ldots,30\}$. First, we note that it is not sensible to estimate them via
\begin{equation}\label{eq:estimates1}
\hat\alpha_t=\frac{R_{t-1}}{N_0}\qquad\text{and}\qquad \hat\beta_t=\frac{I_{t-1}}{N_0}
\end{equation}
(see, e.g., \cite[pages 3--4]{AldibasiEtAl}); indeed, such an estimation gives unreasonably low values of $\cR^{\eff}_t$ throughout the period, as plotted in Figure \ref{fig:preliminaryReff}. The computation, however, also shows that there is no significant difference in the values of $\cR^{\eff}_t$ if the quantities $S_t$ and $N_t$ in the formula \eqref{eq:Reff} are replaced with $\hat S_t$ and $\hat N_t:=\hat S_t+\hat I_t+\hat R_t$. Notice that, while $N_t$ is constant for every $t$, this is not the case for $\hat N_t$.

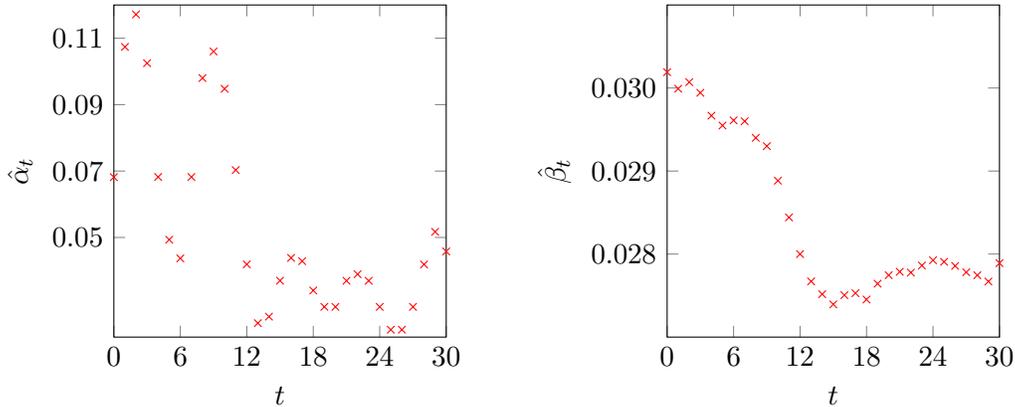
\begin{figure}
\begin{tikzpicture} \pgfkeys{/pgf/number format/.cd,fixed,precision=0}
\begin{axis}[
	xmin=0,
	xmax=30,
	ymin=0.04,
	ymax=0.14,
    xtick={0,6,12,18,24,30},
    ytick={0.07,0.09,0.11,0.13},
    yticklabels={0.05,0.07,0.09,0.11,0.13},
	xlabel={$t$},
	ylabel={$\hat\alpha_t$},
	ylabel near ticks,
	samples=100,
	width=6cm,
	height=6cm
]
\addplot[only marks,mark=x,red] coordinates {(0,0.08816448) (1,0.1273748) (2,0.13716169) (3,0.12248198) (4,0.08821205) (5,0.06931798) (6,6.37e-02) (7,8.82e-02) (8,1.18e-01) (9,1.26e-01) (10,0.1147683) (11,0.09028618) (12,6.19e-02) (13,0.04419672) (14,0.0461684) (15,0.05698091) (16,0.06385045) (17,0.06286678) (18,0.05403032) (19,0.04910888) (20,0.04910371) (21,0.05695789) (22,0.05892289) (23,0.05695362) (24,0.0490947) (25,0.04222215) (26,0.04222408) (27,0.04910187) (28,0.06187049) (29,0.07169682) (30,0.06578914)};
\end{axis}
\end{tikzpicture}\qquad\begin{tikzpicture} \pgfkeys{/pgf/number format/.cd,fixed,precision=0}
\begin{axis}[
	xmin=0,
	xmax=30,
	ymin=0.027,
	ymax=0.031,
    xtick={0,6,12,18,24,30},
    ytick={0.028,0.029,0.030},
    yticklabels={0.028,0.029,0.030},
	xlabel={$t$},
	ylabel={$\hat\beta_t$},
	ylabel near ticks,
	samples=100,
	scaled y ticks=false,
	width=6cm,
	height=6cm
]
\addplot[only marks,mark=x,red] coordinates {(0,0.03019077) (1,0.02999159) (2,0.03007087) (3,0.02994273) (4,0.02966688) (5,0.02954828) (6,0.0296102) (7,2.96e-02) (8,2.94e-02) (9,2.93e-02) (10,0.0288839) (11,0.02844151) (12,2.80e-02) (13,0.02767088) (14,0.02751676) (15,0.02739403) (16,0.02750517) (17,0.02752874) (18,0.02745276) (19,0.02764322) (20,0.02774555) (21,0.02778676) (22,0.0277769) (23,0.02786215) (24,0.02792521) (25,0.0279068) (26,0.02785747) (27,0.02778096) (28,0.02774485) (29,0.02766954) (30,0.0278905)};
\end{axis}
\end{tikzpicture}
\caption{\label{fig:alphabeta} Plots of the points $\left\{\left(t,\hat\alpha_t\right)\right\}_{t=0}^{30}$ (left) and $\left\{\left(t,\hat\beta_t\right)\right\}_{t=0}^{30}$ (right) resulting from the L-BFGS-B algorithm.}
\end{figure}

We therefore apply the following alternative method to estimate the values of $\hat\alpha_t$ and $\hat\beta_t$. On any given day $t\in\{0,\ldots,30\}$, we let the pair $\left(\hat\alpha_t,\hat\beta_t\right)$ be the least-squares solution of the system of 90 linear equations in $\hat\alpha_t$ and $\hat\beta_t$ constructed by substituting into \eqref{eq:discretemodel} the parameter values \eqref{eq:fixedparameters} as well as the values of $\rho_s$ and $\left(\hat S_s,\hat I_s,\hat R_s\right)=\left(S_s,I_s,R_s\right)$ for every $s\in\{t-90,\ldots,t-1\}$. Such a method is referred to as the \textit{L-BFGS-B algorithm} \cite{ByrdEtAl,CoppolaStewart,FeiRongWangWang}. The resulting points $\left\{\left(t,\hat\alpha_t\right)\right\}_{t=0}^{30}$ and $\left\{\left(t,\hat\beta_t\right)\right\}_{t=0}^{30}$ are plotted in Figure \ref{fig:alphabeta}.

Subsequently, computing the values of $\cR^{\eff}_t$ for $t\in\{0,\ldots,30\}$ using the formula \eqref{eq:Reff}, we obtain the points $\left\{\left(t,\cR^{\eff}_t\right)\right\}_{t=0}^{30}$ plotted in Figure \ref{fig:Reff}. It is apparent from the figure that, for $t\in\{0,\ldots,30\}$,
$$\cR^{\eff}_t\in(0.40,0.72)\qquad\text{and}\qquad \rho_t\in\begin{cases}
(0.2,1], &\text{if }t\in\{0,\ldots,17\};\\
[0,0.2], &\text{if }t\in\{18,\ldots,30\}.
\end{cases}$$

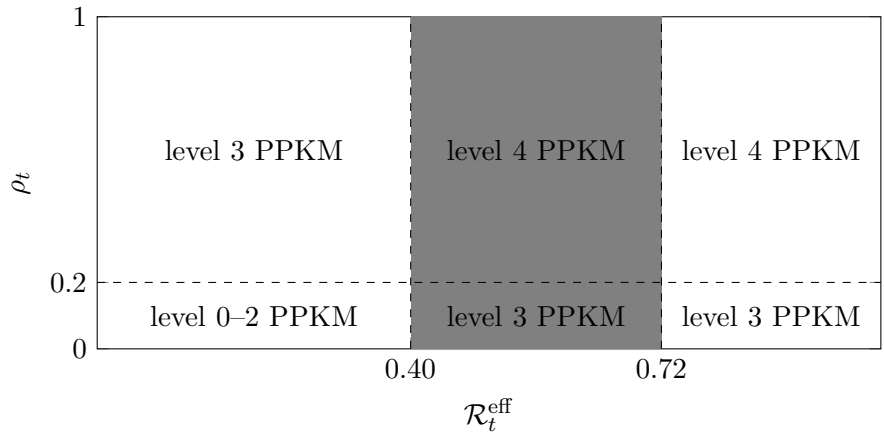
\begin{figure}
\begin{tikzpicture} \pgfkeys{/pgf/number format/.cd,fixed,precision=0}
\begin{axis}[
	xmin=0,
	xmax=1,
	ymin=0,
	ymax=1,
    xtick={0.40,0.72},
    ytick={0,0.2,1},
    xticklabels={0.40,0.72},
    yticklabels={0,0.2,1},
	xlabel={$\cR^{\eff}_t$},
	ylabel={$\rho_t$},
	ylabel near ticks,
	samples=100,
	width=12cm,
	height=6cm
]
\fill[gray] (axis cs:0.4,0) rectangle (axis cs:0.72,1);
\draw[dashed] (axis cs:0,0.2) -- (axis cs:1,0.2);
\draw[dashed] (axis cs:0.4,0) -- (axis cs:0.4,1);
\draw[dashed] (axis cs:0.72,0) -- (axis cs:0.72,1);
\node at (axis cs:0.56,0.1) {level 3 PPKM};
\node at (axis cs:0.56,0.6) {level 4 PPKM};
\node at (axis cs:0.86,0.1) {level 3 PPKM};
\node at (axis cs:0.86,0.6) {level 4 PPKM};
\node at (axis cs:0.2,0.1) {level 0--2 PPKM};
\node at (axis cs:0.2,0.6) {level 3 PPKM};
\end{axis}
\end{tikzpicture}
\caption{\label{fig:design} Our design of governmental policies: the appropriate level(s) of PPKM which should be enforced in the situation described by every point $\left(\cR^{\eff}_t,\rho_t\right)\in[0,\infty)\times[0,1]$.}
\end{figure}

As mentioned in the previous section, $t\in\{0,\ldots,17\}$ is a period of level 4 PPKM, while $t\in\{18,\ldots,30\}$ is a period of level 3 PPKM. We see therefore that the downgrading of the PPKM level from 4 to 3 on day $t=17$ (August 24th) is not due to the values of $\cR^{\eff}_t$, but to those of $\rho_t$ which undergo a consistent decrease. Consequently, let us assign to the two rectangles $(0.40,0.72)\times (0.2,1]$ and $(0.40,0.72)\times [0,0.2]$ on the $\cR^{\eff}\rho$-plane  (shaded in Figure \ref{fig:design}) the policies of level 4 PPKM and level 3 PPKM, respectively.


To complete our policy design, let us assign an appropriate level of PPKM to each of the four unshaded rectangles: $[0,0.40]\times [0,0.2]$, $[0,0.40]\times (0.2,1]$, $[0.72,\infty)\times [0,0.2]$, and $[0.72,\infty)\times (0.2,1]$, as follows.
\begin{enumerate}[topsep=5pt,itemsep=5pt]
\item[(i)] The rectangle $[0.72,\infty)\times (0.2,1]$ represents a pandemic situation which is worse than that represented by $(0.40,0.72)\times (0.2,1]$, to which we have already assigned level 4. Since there is no PPKM level higher than 4, it makes sense to also assign level 4 to $[0.72,\infty)\times (0.2,1]$.
\item[(ii)] Let us now consider the rectangles $[0.72,\infty)\times [0,0.2]$ and $[0,0.40]\times (0.2,1]$. Compared to $[0.72,\infty)\times (0.2,1]$, in the former rectangle, the reproduction number is on the same range while the bed-occupancy rate is lower, whereas in the latter, the bed-occupancy rate is on the same range while the reproduction number is lower. We thus assign to each of these rectangles the lower level 3. For the former rectangle, any lower level is not recommended; this is to avoid any drastic increase of the bed-occupancy rate due to the high level of reproduction number. By contrast, for the latter, a lower level could be considered, since the currently high bed-occupancy rate can be expected to decline due to the low level of reproduction number.
\item[(iii)] Finally, in the rectangle $[0,0.40]\times [0,0.2]$, where the reproduction number and the bed-occupancy rate are both low, it is fitting to assign level 2 or any lower level of PPKM.
\end{enumerate}
Therefore, in a period where level 3 PPKM is being enforced, the PPKM level could be downgraded to 2 only once \textit{both} the reproduction number and the bed-occupancy rate have reached notably low values (under 0.40 and under 0.2, respectively), and should be upgraded to 4 if \textit{both} of them have exceeded the same respective thresholds (over 0.40 and over 0.2, respectively).

\section{Conclusions and future research}\label{sec:conclusions}

We have constructed a design of governmental policies for the eradication of COVID-19 in DKI Jakarta, consisting in a quantitative method to determine the appropriate PPKM level on any given day $t$, which can be summarised as follows. First, we obtain the present hospitals' bed-occupancy rate $\rho_t$ and effective reproduction number $\cR^{\eff}_t$. Then, we check the validity of two conditions: (i) $\rho_t>0.2$ and (ii) $\cR^{\eff}_t>0.4$. A level 4 PPKM should be enforced if and only if both (i) and (ii) are true, a level 3 PPKM should be enforced if exactly one of (i) and (ii) is true, and a lower level of PPKM can be enforced if both (i) and (ii) are false.

This research can be developed in the following ways. First, one could carry out similar analysis using more advanced mathematical models, taking into account additional aspects such as the possibility of reinfection \cite[section 4]{YongOwenHoseana} or the number of vaccinated individuals, the latter being expected to play a significant role in the determination of the PPKM level \cite{Office}. Improvements could also be made on our parameter estimation methods. Indeed, our estimated values of recovery and incidence rates on any given day could be compared to those estimated using, rather than the data of the previous 90 days, e.g., those of both the previous and following weeks. Finally, our method can also be applied to provinces other than DKI Jakarta.


\end{document}